\documentclass[11pt]{article}

\usepackage{amsmath}
\usepackage{amssymb}

 \newcommand{\beq}{\begin{equation}}
\newcommand{\eeq}{\end{equation}}

\newtheorem{theorem}{Theorem}[section]

\newtheorem{lemma}[theorem]{Lemma}
\newtheorem{corollary}[theorem]{Corollary}

\newcommand{\Rn}{{ \mathbb R}^n}

\newcommand{\Cn}{{\mathbb  C\sp n}}

\newcommand{\D}{{\mathbb D}}

\newcommand{\F}{{\mathcal F}}
\newcommand{\E}{{\mathcal E}}

\newcommand{\cO}{{\mathcal O}}

\newcommand{\cI}{{\cal I}}

\newcommand{\m}{{\mathfrak m}}

\newcommand{\PSH}{{\operatorname{PSH}}}

\newcommand{\codim}{{\operatorname{codim}}}

\begin{document}


\begin{center}
{\Large\bf Extremal cases for the log canonical threshold}
\end{center}

\begin{center}
{\large Alexander Rashkovskii}
\end{center}

\vskip1cm

\begin{abstract}
We show that a recent result of Demailly and Pham Hoang Hiep \cite{DH} implies a description of plurisubharmonic functions with given Monge-Amp\`ere mass and smallest possible log canonical threshold. We also study an equality case for the inequality from \cite{DH}.
\end{abstract}

\section{Introduction and statement of results}

Let $\PSH_0$ denote the collection of of germs of all
functions plurisubharmonic at the origin of $\Cn$.
A basic characteristic of singularity of $u\in
\PSH_0$ is its {\it Lelong number}
$$ \nu_u=\nu_u(0)=\liminf{u(z)}/{\log|z|},\   z\to0.$$
One more characteristic, introduced in various contexts by several authors (first, probably, in \cite{Sk}) and attracted recently considerable attention (e.g., \cite{ACKPZ}, \cite{Be}, \cite{D8}, \cite{DK}, \cite{DH}, \cite{FaJ2}, \cite{H2}, \cite{Kis3}), is the {\it integrability index} (at $0$)
$$ \lambda_u= \inf\{\lambda>0:e^{-u/\lambda}\in L^2_{loc}(0)\}.$$
For an ideal $\cI=\cI(f_1,\ldots,f_m)\subset\cO_0$ generated by analytic germs $f_1,\ldots,f_m$, the value $c({\cI})=\lambda_{\log|f|}^{-1}$ is the {log canonical threshold} of  ${\cI}$.
Accordingly, $c_u=\lambda_{u}^{-1}$ is called the {\it log canonical threshold} of $u$.

A classical result due to Skoda \cite{Sk} states that
\begin{equation}\label{skoda}
\nu_u^{-1}\le c_u\le n\,\nu_u^{-1},
\end{equation}
the extremal situations (equalities) being realized, for example, for
$u=\log|z_1|$ (for the first inequality) and $u=\log|z|$ (for the second one). A description of all functions $u$ with $c_u= n\nu_u^{-1}$ was given in \cite{R10}. The other extremal relation seems to be more involved. The only known to us result in this direction concerns the case $n=2$, where the functions satisfying $c_u=\nu_u^{-1}$ are proved in \cite{FaJ2} to be of the form $u=c\log|f|+v$, where $f$ is an analytic function regular at $0$, and $v\in\PSH_0$ has zero Lelong number at $0$.

In this note, we concentrate on lower bounds for the log canonical threshold, with the main focus when the inequalities become equalities.

In  \cite{dFEM} and \cite{M}, the log canonical threshold of a zero dimensional ideal ${\mathcal I}\subset\cO_0$ was related to its Samuel multiplicity $e({\cI})$:
\begin{equation}\label{dFEM} c({\cI})\ge {n}{e(\cI)^{-1/n}},\end{equation}
with an equality if and only if the integral closure of $\cI$ is a power of the maximal ideal $\m_0\subset\cO_0$.
It was used by Demailly \cite{D8} for a corresponding bound for plurisubharmonic functions $u$ with isolated singularity at $0$, and extended then by Zeriahi \cite{Ze} to all $u$ with $(dd^cu)^n$ well defined (more precisely, for all $u$ from the Cegrell class $\E$ \cite{Ce04}),
\begin{equation}\label{dem} c_u\ge n\,{e_n(u)^{-1/n}}.\end{equation}
Here $e_k(u)$ are the Lelong numbers of the currents  $(dd^c u)^k$ at $0$:
$$ e_k(u)=(dd^c u)^k\wedge(dd^c\log|z|)^{n-k}(0), \quad 1\le k\le n,$$
and $d=\partial + \bar\partial$, $d^c= (\partial -\bar\partial)/2\pi i$.
 The Cegrell class $\F(D)$ is formed by limits of decreasing sequences of bounded plurisubharmonic functions $u_j$ in $D$ such that $u_j=0$ on $\partial D$ and $\sup_j\int_D (dd^cu_j)^n<\infty$, and  $u\in\E(D)$ if for any $K\Subset D$ one can find $v\in\F(D)$ such that $u=v$ on $K$, see  \cite{Ce04}. In particular, all negative plurisubharmonic functions that are bounded outside a compact subset of $D$,  belong to $\E(D)$.

Note that $e_1(u)=\nu_u$. When $\cI=\cI(f_1,\ldots,f_p)\subset\cO_0$ is a zero dimensional ideal, then $e_n(\log|f|)=e(\cI)$, see \cite{D8}. If $\codim\, V(\cI)=k$, the values $e_k(\log|f|)$ are {mixed Rees multiplicities} $e_k(\cI)$ of $\cI$ and the maximal ideal $\mathfrak{m}_0$  considered, e.g., in \cite{BA}.

A direct proof of Demailly's inequality (\ref{dem}) without using (\ref{dFEM}) was obtained in \cite{ACKPZ}. In \cite{D8}, the question of equality in (\ref{dem}) has been raised, and it was conjectured that, similarly to the analytic case $u=\log|f|$, the extremal functions should be plurisubharmonic functions with logarithmic singularity at $0$.

In \cite{R10}, Demailly's inequality was used to get the `intermediate' bounds
\beq\label{eq:rel} c_u\ge k\,e_k(u)^{-1/k},\quad 1\le k\le l,\eeq
where $l$ is the codimension of an analytic set $A$ such that $u^{-1}(-\infty)\subset A$. None of the bounds for different values of $k$ can be deduced from the others.

In a recent paper \cite{DH}, an optimal bound for the integrability index in terms of the Lelong numbers was obtained:
if $u\in\E$ near $0$ and $e_1(u)>0$, then
\beq\label{DH} c_u\ge E_n(u):=\sum_{1\le k\le n}\frac{e_{k-1}(u)}{e_{k}(u)},\eeq
where $e_0(u)=1$.
It is easy to see that this bound implies all the relations (\ref{eq:rel}) for the case of $l=n$ (that is, for $u$ with isolated singularity).
Here we will show that it also gives an answer to the aforementioned question on equality in (\ref{dem}).

To state it, we need the following notion from \cite{R7}. Let $D$ be a bounded, hypeconvex neighborhood of $0$. Given a function $u\in\PSH^-(D)$ (negative and plurisubharmonic in $D$), its {\it greenification} $g_u$ at $0$ is the regularized upper envelope of all functions $v\in\PSH^-(D)$ such that $v\le u+O(1)$ near $0$.

The greenification of $\log|z|$ is the standard pluricomplex Green function with pole at $0$. For $u$ satisfying $(dd^cu)^n=0$ on a punctured neighborhood of the origin, $g_u$ is the Green function in the sense of Zahariuta \cite{Za0}. The greenification of a {\it multi-circled singularity} $u(z)=u(|z_1|,\ldots,|z_n|)+O(1)$ in the unit polydisk $\D^n$ is the so-called {\it indicator}: a multi-circled function satisfying $g_u(|z_1|^c,\ldots,|z_n|^c)=c\,g_u(z)$ $\forall c>0$ \cite{R10}.

One has always $(dd^cg_u)^n=0$ on $\{g_u>-\infty\}$. Evidently, $ g_u\ge u$, while the relation $g_u=u+O(1)$ need not be true. Nevertheless, the greenification keeps the considered characteristics of singularity:

\medskip
\begin{lemma}\label{lem:gr}
  Let $u\in\PSH_0$ and let $g_u$ be its greenification on a bounded hyperconvex neighborhood $D$ of $0$. Then $\lambda_{g_u}=\lambda_u$.  If, in addition, $u\in \E$ on a neighborhood of $0$, then
$g_u\in\F(D)$, $(dd^c g_u)^n=0$ on $D\setminus\{0\}$, and $e_k(g_u)=e_k(u)$ for all $k$.
\end{lemma}
\medskip

Therefore, the only information on asymptotic behavior of $u$ one can expect from the values of $c_u$ and $e_k$ is the one on its greenifications $g_u$.

\medskip
\begin{theorem}\label{theo1} For any $u\in\E$ near $0$, the relation
$c_u= n\,{e_n(u)^{-1/n}}$ holds if and only if its greenification for some (and then for any) bounded hyperconvex domain $D$ satisfies $g_u=e_1(u)\log|z|+O(1)$ as $z\to0$.
\end{theorem}
\medskip

\begin{corollary} Let $u\in\F(D)$, $e_1(u)=1$, and
$ \int_D (dd^cu)^n=(n\lambda_u)^n$. Then
$u$ is the pluricomplex Green function for $D$ with logarithmic singularity at $0$.
\end{corollary}
\medskip

In the case of analytic singularities $u=\log|f|$, statement (i) of Theorem~\ref{theo1} recovers the aforementioned result from \cite{dFEM} on equality in the bound for log canonical thresholds.

\medskip
Next question is when equalities in (\ref{eq:rel}) and (\ref{DH}) occur. Moreover, the latter bound can be extended to the case of functions not from $\E$, which rises a question on the equality cases.

\medskip
\begin{theorem}\label{theo2}
If $u\in\PSH_0$ is locally bounded outside an analytic set of codimension $l>1$, then
\beq\label{DH1} c_u\ge E_l(u):=\sum_{1\le k\le l}\frac{e_{k-1}(u)}{e_{k}(u)}.\eeq
\end{theorem}
\medskip

(Note that relation (\ref{DH1}) for $l=1$ is the lower bound in Skoda's inequalities (\ref{skoda}) and it does not require any assumption on $u$.)

\medskip

For multi-circled singularities $\varphi(z)=\varphi(|z_1|,\ldots,|z_n|)+O(1)$ and any $l$, it was proved in \cite{R10} that the relation $c_\varphi= l\,e_l(\varphi)^{-1/l}$ holds if and only if its greenification $g_\varphi$ in  $\D^n$ equals $e_1(z)\max_{j\in J}\log|z_j|$ for an $l$-tuple $J\subset\{1,\ldots,n\}$.

\medskip

\begin{theorem}\label{mult_DH} If a multi-circled plurisubharmonic singularity $\varphi$ satisfies
$c_\varphi= E_l(\varphi)$,
then
\beq\label{simplex} g_\varphi(z)=\max_{j\in J}\frac{\log|z_{j}|}{a_j}\eeq
for some $l$-tuple $ J=(j_1,\ldots,j_l)\subset\{1,\ldots,n\}$ and $a_j>0$.
\end{theorem}

\medskip

A characterization of functions  of the form (\ref{simplex}) is that they generate {\it monomial valuations} ${\frak v}_\varphi$ on plurisubharmonic singularities $u$ by ${\frak v}_\varphi(u)=\liminf u(z)/\varphi(z)$ as ${z\to 0}$.
One could ask if the statement of Theorem~\ref{mult_DH} remains true for $\varphi$ generating {\it quasi-monomial valuations}, i.e., monomial ones on birational models \cite{BFaJ}. As the following example shows, the answer is no.

\medskip
{\it Example 2.}
As follows from \cite{FaJ}, the function $\varphi=\log(|z_1^4|+|z_1^3-z_2^2|)$ generates a quasi-monomial valuation. Since $u=\log|z_1^3-z_2^2|\le \varphi\le v=\log(|z_1^4|+|z_1^3|+|z_2^2|)$
and $c_u=c_v= 5/6$, we have $c_\varphi=5/6>E_2(\varphi)=3/4$.

\section{Proofs}

\noindent
1. {\it  Proof of Lemma~\ref{lem:gr}.} Evidently, $c_{g_u}\ge c_u$. By the Choquet lemma, there exists a sequence $u_j$ increasing a.e. to $g_u$ and such that $u_j\le u+O(1)$ and so, $c_{u_j}\le c_{u}$. Semicontinuity theorem \cite{DK} shows then $c_{g_u}\le c_u$.

Let $u\in \E(\omega)$,  $0\in\omega\subset D$. Then there exists $v\in\F(\omega)$ such that $v=u$ near $0$. Furthermore, there exists $w\in\F(D)$ such that $w\le v$ on $\omega$ \cite{CeZe}. Since $w\le g_u$, the function $g_u$ belongs to $\F(D)$. The relation $(dd^c g_u)^n=0$ outside $0$ follows by standard arguments, because maximality of $v\in\E$ on an open set $U$ is equivalent to $(dd^c v)^n(U)=0$.

To prove $e_k(g_u)= e_k(u)$, we take again a sequence $u_j$ increasing a.e. to $g_u$; $u_j$ can be chosen to be from the class $\F(D)$, for otherwise we replace them by $\max\{u_j,w\}$. Therefore, the currents $(dd^cu_j)^k$ converge to $(dd^cg_u)^k$ \cite{Ce12} (the result is stated there only on the convergence of $(dd^cu_j)^n$, while the proof uses induction in the degree $k$). By the semicontinuity theorem for the Lelong numbers \cite{D1}, this implies $\limsup_{j\to\infty} e_k(u_j)\le e_k(g_u)$. On the other hand, the relations $u_j\le u+O(1)\le g_u$ give us, by the comparison theorem for the Lelong numbers \cite{D1}, $e_k(u)\le \limsup e_k(u_j)$ and $e_k(g_u)\le e_k(u)$.
%
{\hfill$\square$\rm \vskip 2mm}

\medskip\noindent
2. Further proofs are based essentially on estimate (\ref{DH}) and the following uniqueness result.

\medskip
 \begin{lemma}\label{lem:uniq} {\rm (\cite[Thm. 3.7]{ACCP}; for greenifications of isolated singularities, \cite[Lem. 6.3]{R7})} If $u,v\in\F(D)$ are such that $u\le v$ and $(dd^cu)^n=(dd^cv)^n$, then $u=v$. As a consequence, if $u,v\in\E$,  $u\le v +O(1)$ near $0$, and $e_n(u)=e_n(v)$, then $g_u=g_v$.
\end{lemma}

\medskip
\noindent
3. {\it  Proof of Theorem~\ref{theo1}.} By Lemma \ref{lem:gr}, we can assume $u=g_u$. Relation (\ref{DH}) gives us
$E_n(u)=n\,{e_n(u)^{-1/n}}$,
and by the arithmetic-geometric mean theorem we get then
$$\frac{e_{k-1}(u)}{e_{k}(u)}=\frac{e_{l-1}(u)}{e_{l}(u)}$$
for any $k,l\le n$, which implies $e_n(u)=[e_1(u)]^n$.
Let $v=e_1(u)G$, where $G$ denotes the pluricomplex Green function for $D$ with logarithmic pole at $0$. Since  $u\in \PSH^-(D)$ satisfies $u\le e_1(u)\log|z|+O(1)$ as $z\to0$, we have $u\le v$ on $D$, while $e_n(u)=e_n(v)$. By Lemma~\ref{lem:uniq} we conclude then $u=v$.
{\hfill$\square$\rm \vskip 2mm}

\medskip
\noindent
4. {\it  Proof of Theorem~\ref{theo2}.} The restriction $u_L$ of $u$ to a generic $l$-dimensional subspace $L\in G(l,n)$ has isolated singularity at $0$ and, by Siu's theorem, $e_k(u_L)=e_k(u)$. By \cite[Prop.~2.2]{DK}, we have also $c_u\ge c_{u_L}$. Therefore, we can apply (\ref{DH}) to $u_L$ and get the bound (\ref{DH1}).
{\hfill$\square$\rm \vskip 2mm}

\medskip\noindent
5. {\it Proof of Proposition~\ref{mult_DH}.}
By considering again restriction to a generic $l$-dimensional coordinate plane, we can assume $l=n$ and $\varphi$ to coincide with its greenification $g_\varphi$ in $\D^n$.

As was proved in \cite{DH}, the bound (\ref{DH}) for multi-circled functions follows from the inequality
$$\varphi\le \Phi(z):=|\varphi(z^*)|\max_{j}\frac{\log|z_{j}|}{a_j}$$ with $a_j=|\log|z_j^*||$ (e.g., \cite[Prop. 3]{R}), where $z^*\in\Pi=\{z: \: |z_1\cdot\ldots\cdot z_n|=1/e\}$ is chosen such that
$|\varphi(z^*)|=|\min\{\varphi(z): z\in\Pi \}|=\lambda_\varphi$ \cite[Thm. 5.8]{Kis3}. Namely,  $E_k(\varphi)\le E_k(\Phi)$ for all $k$ and $c_\varphi=c_\Phi=E_n(\Phi)$. Therefore, $c_\varphi=E_n(\varphi)$ implies $E_n(\varphi)=E_n(\Phi)$.

Following \cite{DH}, we set $t_0=1$ and consider the function $f(t)=\sum_1^n t_{j-1}/t_j$
on the convex set $\{t\in\Rn_+:\: t_j^2\le t_{j-1}t_{j+1}\}$. The function is decreasing in each variable $t_j$ and strictly decreasing in $t_n$.
Note that $E_n(v)=f_n(e_1(v),\ldots,e_n(v))$ for any $v\in\PSH_0$ with isolated singularity. If $e_n(\varphi)>e_n(\Phi)$, then we would have $E_n(\varphi)<E_n(\Phi)$, which is not true. Therefore, $e_n(\varphi)=e_n(\Phi)$ and, since $\varphi\le\Phi$, we get $\varphi=\Phi$ by Lemma~\ref{lem:uniq}.

\end{document}